# SOLITARY WAVES FOR MAXWELL-SCHRÖDINGER EQUATIONS

GIUSEPPE MARIA COCLITE AND VLADIMIR GEORGIEV

**Abstract**

In this paper we study the solitary waves for the coupled Schrödinger - Maxwell equations in three-dimensional space. We prove the existence of a sequence of radial solitary waves for these equations with a fixed $L^2$ norm. We study the asymptotic behavior and the smoothness of these solutions. We show also the fact that the eigenvalues are negative and the first one is isolated.

## 1. Introduction

The wave function $\psi(t, x), t \in \mathbb{R}, x \in \mathbb{R}^3$ of a non - relativistic charged particle (say one electron) is a solution of the classical Schrödinger equation

$$i\hbar \partial_t \psi + \mathbb{H}\psi = 0, \tag{1.1}$$

where $\hbar$ is the Planck constant and $\mathbb{H}$ is the Hamiltonian, i.e. (see [15], 2.17, p.20)

$$\mathbb{H} = -\frac{\hbar^2}{2m}\Delta + \mathbb{W}. \tag{1.2}$$

Here $m > 0$ is the mass of the particle and $\mathbb{W}$ is the operator of the total - potential energy. In the case of one electron with charge $e > 0$ and one nucleus with charge number $z$ (i.e. $z$ is the number of the protons in the nucleus) we have (see (1.4), p.14 in [15])

$$\mathbb{W} = -e^2 z |x|^{-1} + \mathbb{B}, \tag{1.3}$$

where $-e^2 z |x|^{-1}$ is the electron - nucleus attractive potential energy, while $\mathbb{B}$ is the electron - electron repulsive energy. As in [15] we choose units in which

$$\frac{\hbar^2}{2m} = 1, \quad e = 1.$$

The interaction between the electromagnetic field and the wave function related to a quantistic non-relativistic charged particle (say one electron) is described by the Maxwell - Schrödinger system. More precisely, let $\psi = \psi(x, t)$ be the wave function and let $\mathcal{A} = (A^0, A^1, A^2, A^3)$ be the electromagnetic potentials of a charged

The authors are partially supported by Research Training Network (RTN) HYKE, financed by the European Union, contract number : HPRN-CT-2002-00282.

non- relativistic particle. The the corresponding Maxwell - Schrödinger system (in Lorenz gauge) have the form

$$\partial_{tt}\mathcal{A} - \Delta\mathcal{A} = \mathcal{J},$$

$$i\hbar\partial_{t,\mathcal{A}}\psi - \frac{1}{2}\Delta_{\mathcal{A}}\psi - V(x)\psi = 0,$$

(1.4) $$\partial_t A^0 + \sum_{k=1}^{3} \partial_{x_k} A^k = 0,$$

where

(1.5)
$$\partial_{t,\mathcal{A}} = \partial_t + iA_0, \ \Delta_{\mathcal{A}} = \sum_{k=1}^{3} \partial^2_{k,\mathcal{A}},$$
$$\partial_{k,\mathcal{A}} = \partial_{x_k} + iA_k,$$
$$\mathcal{J} = (J_0, J_1, J_2, J_3),$$
$$J_0 = 4\pi|\psi|^2, J_k = 4\pi\mathrm{Im}\left(\bar{\psi}\partial_{k,\mathcal{A}}\psi\right).$$

Here $V = V(r)$ is the (external) potential of the electron - nucleus attractive potential energy. From (1.3) we have

(1.6) $$V(|x|) = \frac{z}{|x|},$$

where $z$ denotes the number of protons in a nucleus ($z = 1$ for hydrogen).

We consider the electrostatic case, namely

$$A_0 = \varphi(x)/\hbar, \qquad A_j(x) = 0, \qquad j = 1,2,3, \qquad x \in \mathbb{R}^3,$$

and for standing wave function

$$\psi(x,t) = u(x)e^{i\omega t/\hbar}, \qquad x \in \mathbb{R}^3, \ t \in \mathbb{R},$$

where $\omega \in \mathbb{R}$ and $u$ is real valued. In the case of zero magnetic field (i.e. $A_j = 0, j = 1,2,3,$ ) and non - trivial electric potential $\varphi(x)$ the Maxwell - Schrödinger system has the form (see [6])

$$-\frac{1}{2}\Delta u - \varphi u - V(x)u = \omega u, x \in \mathbb{R}^3,$$
$$\Delta\varphi = 4\pi u^2, \ x \in \mathbb{R}^3,$$
(1.7) $$\int_{\mathbb{R}^3} u^2 = N,$$

where the last equation is due to the probabilistic interpretation of the wave function. The number $N$ of electrons shall be assumed less or equal the number of protons, i.e.

(1.8) $$N \leq z.$$

The equations in (1.7) have a variational structure, in fact they are the Euler - Lagrange equations related to the functional:

(1.9) $$F(u,\varphi) = \frac{1}{4}\int_{\mathbb{R}^3}|\nabla u|^2 dx - \frac{1}{2}\int_{\mathbb{R}^3}\varphi u^2 dx - \frac{1}{2}\int_{\mathbb{R}^3}V(x)u^2 dx - \frac{1}{16\pi}\int_{\mathbb{R}^3}|\nabla\varphi|^2 dx.$$

It is easy to see that this functional is well - defined, when

$$u \in H^1(\mathbb{R}^3), \quad \int_{\mathbb{R}^3}|\nabla\varphi|^2 dx < \infty.$$

This functional is strongly indefinite, which means that $F$ is neither bounded from below nor from above and this indefiniteness cannot be removed by a compact perturbation. Moreover $F$ is not even. Later on ( see (2.5) ) we shall introduce a functional $J(u)$ that is bounded from below and such that the critical points of $J$ can be associated with the critical points of $F(u, \varphi)$.

The first natural question is connected with the simplest case $V \equiv 0$ (that is $z = 0$), namely

$$-\frac{1}{2}\Delta u - \varphi u = \omega u, x \in \mathbb{R}^3,$$
(1.10)
$$\Delta \varphi = 4\pi u^2, \ x \in \mathbb{R}^3.$$

It is well-known that the similar physical model of Maxwell - Dirac system with zero external field admits solitary solutions (see [13]), i.e. nontrivial solutions in the Schwartz class $S(\mathbb{R}^3)$.

Our first result is the following.

**Theorem 1.1.** *Let $(u, \varphi, \omega)$ be a solution of (1.10) such that $u$, $\varphi$ radial and*

$$u \in H^1(\mathbb{R}^3), \quad \int_{\mathbb{R}^3} |\nabla \varphi|^2 dx < \infty.$$

*Then*

$$u \equiv \varphi \equiv 0.$$

The above result shows that the Schrödinger - Maxwell equations with zero potential have only the trivial solution. This fact justifies the study of the Schrödinger - Maxwell equations with nonzero external potential.

We shall look for soliton type solutions $u$, i.e. very regular solutions decaying rapidly at infinity. First, we establish the existence of $H^1$ radially symmetric solutions.

**Theorem 1.2.** *Under the assumptions (1.6) and (1.8), there exists a sequence of real negative numbers $\{\omega_k\}_{k \in \mathbb{N}}$ so that*

$$\omega_k \longrightarrow 0$$

*and for any $\omega_k$ there exists a couple $(u_k, \varphi_k)$ of solutions of (1.7) such that*

$$u_k \in H^1(\mathbb{R}^3), \quad \int_{\mathbb{R}^3} |\nabla \varphi_k|^2 dx < \infty.$$

*Moreover $u_k, \varphi_k$ are radially symmetric functions.*

A more precise information about the localization of the eigenvalues $\omega$ is given in the following.

**Theorem 1.3.** *Assume (1.6) and (1.8). Let $(u, \varphi, \omega)$ be a nontrivial solution of the equations in (1.7) such that $u$, $\varphi$ radial and*

$$u \in H^1(\mathbb{R}^3), \quad \int_{\mathbb{R}^3} |\nabla \varphi|^2 dx < \infty.$$

*Then we have*

(1.11)
$$\omega < 0.$$

On the other hand, the solutions constructed in Theorem 1.2 are only radial ones. Therefore, it remains as an open problem the existence of non-radial solutions.

The qualitative properties of the solutions are described in the following.

**Theorem 1.4.** *Under the assumptions (1.6), if $(u, \varphi, \omega)$ is a solution of (1.7) with $u$ and $\varphi$ radially symmetric maps and such that*

$$u \in H^1(\mathbb{R}^3), \quad \int_{\mathbb{R}^3} |\nabla \varphi|^2 dx < \infty,$$

*then*

*a) $u(r) \in C^\infty([0, 1]), \varphi(r) \in C^\infty([0, 1])$;*

*b) if $N = z$ then $u \in S(|x| > 1)$, with $S(|x| > 1)$ being the Schwartz class of rapidly decreasing functions.*

*Remark* 1.1. The property b) in the above theorem shows that the soliton behavior of the solutions can be established, when the neutrality condition $N = z$ is satisfied. The physical importance of the neutrality condition is discussed in [15] ( see (5.2) page 24 in [15]).

Finally the topological properties of the set of the solutions are stated in the following.

**Theorem 1.5.** *Under the assumptions (1.6) and (1.8), let $(u, \varphi, \omega)$ be a solution of (1.7) such that $\omega < 0$ is the first eigenvalue, $u$ and $\varphi$ are radially symmetric maps and such that*

$$u \in H^1(\mathbb{R}^3), \quad \int_{\mathbb{R}^3} |\nabla \varphi|^2 dx < \infty.$$

*Then the solution $(u, \varphi, \omega)$ is isolated, i.e. there exists a neighborhood $U \subset H^1(\mathbb{R}^3)$ of $u$, one $W$ of $\varphi$ such that*

$$\int_{\mathbb{R}^3} |\nabla \phi|^2 dx < \infty, \qquad \phi \in W$$

*and one $\Omega \subset \mathbb{R}$ of $\omega$ such that each $(v, \phi, \lambda) \in U \times W \times \Omega$ with $(v, \phi, \lambda) \neq (u, \varphi, \omega)$, $v$ and $\phi$ radially symmetric maps satisfying the following*

$$\int_{\mathbb{R}^3} |v|^2 dx = N,$$

*is not a solution of (1.7).*

For the sake of completeness we want to recall that the existence of solitary waves has been studied by Benci and Fortunato (see [6]) in the case in which the charged particle "lives" in a bounded space region $\Omega$.

Moreover, the Maxwell equations coupled with nonlinear Klein-Gordon equation, with Dirac equation, with nonlinear Schrödinger equation and with the Schrödinger equation under the action of some external potential have been studied respectively in [7, 13, 9, 10, 11].

Finally we recall the classical papers [4, 5, 12]

The plan of the work is the following. In Section 2 we prove some preliminary variational results, that permit to reduce (1.7) to a single equation. Moreover we show the variational structure of the problem. In Section 3 we prove some topological properties of the energy functional associated to (2.4). In Section 4 we prove Theorem 1.1 and 1.3. In Section 5, 6 and 7 we prove Theorem 1.2, 1.4 and 1.5, respectively.

## 2. The Variational Setting

In this section we shall prove a variational principle that permits to reduce (1.7) to the study of the critical points of an even functional, which is not strongly indefinite. To this end we need some technical preliminaries.

We define the space $\mathcal{D}^{1,2}(\mathbb{R}^3)$ as the closure of $C_0^\infty(\mathbb{R}^3)$ with respect to the norm

$$\|u\|_{\mathcal{D}^{1,2}} \doteq \Big(\int_{\mathbb{R}^3} |\nabla u|^2 dx\Big)^{1/2}.$$

The Sobolev - Hardy inequality (see [19]) implies the following.

**Lemma 2.1.** *For all $\rho \in L^{6/5}(\mathbb{R}^3)$ there exists only one $\varphi \in \mathcal{D}^{1,2}(\mathbb{R}^3)$ such that $\Delta \varphi = \rho$. Moreover there results*

(2.1) $$\|\varphi\|_{\mathcal{D}^{1,2}}^2 \leq c\|\rho\|_{L^{6/5}}^2$$

*and the map*

$$\rho \in L^{6/5}(\mathbb{R}^3) \longmapsto \varphi = \Delta^{-1}(\rho) \in \mathcal{D}^{1,2}(\mathbb{R}^3)$$

*is continuous.*

Moreover, the classical Sobolev embedding and a duality argument guarantee the properties

(2.2) $$H^1(\mathbb{R}^3) \subseteq L^q(\mathbb{R}^3) \quad \text{for} \quad 2 \leq q \leq 6$$

$$L^{q'}(\mathbb{R}^3) \subseteq \big(H^1(\mathbb{R}^3)\big)' \quad \text{for} \quad \frac{6}{5} \leq q' \leq 2.$$

Denoting by $H_r^1(\mathbb{R}^3)$ the set of all $H^1$ radial functions. Then the classical Strauss Lemma shows that (see [20] or [8, Theorem A.I'])

(2.3) $\quad H_r^1(\mathbb{R}^3)$ is compactly embedded into $L^q(\mathbb{R}^3), 2 < q < 6.$

By Lemma 2.1 and by using the Sobolev inequalities, for any given $u \in H^1(\mathbb{R}^3)$ the second equation of (1.7) has the unique solution

$$\varphi = 4\pi \Delta^{-1} u^2 \ \big(\in \mathcal{D}^{1,2}(\mathbb{R}^3)\big).$$

For this reason we can reduce the system (1.7) to the equation

(2.4) $$-\frac{1}{2}\Delta u - 4\pi(\Delta^{-1}u^2)u - V(x)u = \omega u, \qquad \int_{\mathbb{R}^3} |u|^2 dx = N.$$

Observe that (2.4) is the Euler-Lagrange equation of the functional

(2.5) $$J(u) = \frac{1}{4}\int_{\mathbb{R}^3} |\nabla u|^2 dx + \pi \int_{\mathbb{R}^3} |\nabla \Delta^{-1} u^2|^2 dx - \frac{1}{2}\int_{\mathbb{R}^3} V(x)|u|^2 dx,$$

constrained on the manifold

$$B \doteq \big\{u \in H^1(\mathbb{R}^3)\big| \|u\|_{L^2}^2 = N\big\}.$$

Note that the functional $J(u)$ can be defined for complex valued $u$, while its critical points are only real-valued.

Given any integer $k \geq 1$ we set

$$H_r^k(\mathbb{R}^3) \doteq \{u \in H^k(\mathbb{R}^3) \mid u(x) = u(|x|), \ x \in \mathbb{R}^3\}.$$

**Lemma 2.2.** *There results:*
  i) *$J$ is even;*

   ii) $J$ is $C^1$ on $H^1(\mathbb{R}^3)$ and its critical points constrained on $B$ are the solutions of (2.4);
   iii) any critical point of $J\big|_{H_r^1(\mathbb{R}^3)\cap B}$ is also a critical point of $J\big|_B$.

PROOF. The proof of i) is trivial.
   Since
$$\frac{d}{d\lambda}\Big(\int_{\mathbb{R}^3}|\nabla\Delta^{-1}|u+\lambda v|^2|^2 dx\Big)\Big|_{\lambda=0} = -4\,\mathrm{Re}\int_{\mathbb{R}^3}(\Delta^{-1}u|v)dx,$$
ii) holds true.

   Now we prove iii). Consider the $O(3)$ group action $T_g$ on $H^1(\mathbb{R}^3)$ defined by
$$T_g u(x) = u(g(x)),$$
where $g \in O(3)$ and $u \in H^1(\mathbb{R}^3)$. Then the conclusion follows by well known arguments (see for example [20]) because $J$ is invariant under the $T_g$ action, namely
$$J(T_g u) = J(u),$$
where $g \in O(3)$ and $u \in H^1(\mathbb{R}^3)$. So, by [17] or [22, Theorem 1.28], iii) is proved. □

3. TOPOLOGICAL RESULTS

In this section we shall prove some topological properties of the functional $J$.

**Lemma 3.1.** *The functional $J$ is weakly lower semicontinuous on $H_r^1(\mathbb{R}^3)$. In particular, the operator*
$$T : u \in H_r^1(\mathbb{R}^3) \longmapsto (\Delta^{-1}u^2)u \in \big(H_r^1(\mathbb{R}^3)\big)'$$
*is compact and the functionals*
$$J_1 : u \in H_r^1(\mathbb{R}^3) \longmapsto \int_{\mathbb{R}^3}|\nabla\Delta^{-1}u^2|^2 dx,$$
$$J_2 : u \in H^1(\mathbb{R}^3) \longmapsto \int_{\mathbb{R}^3}V(x)u^2 dx$$
*are weakly continuous.*

PROOF. We prove that $T$ is compact. Let $\{u_k\} \subset H_r^1(\mathbb{R}^3)$ be bounded. Passing to a subsequence, there exists $u \in H_r^1(\mathbb{R}^3)$ such that
$$u_k \rightharpoonup u \qquad \text{weakly in } H_r^1(\mathbb{R}^3).$$
By (2.1) and Sobolev inequalities (2.2) we see that $\{\Delta^{-1}u_k^2\}$ is bounded in $\mathcal{D}^{1,2}(\mathbb{R}^3)$. Passing to a subsequence, there exists $h \in \mathcal{D}^{1,2}(\mathbb{R}^3)$ such that
(3.1) $$\Delta^{-1}u_k^2 \rightharpoonup h \qquad \text{weakly in } \mathcal{D}^{1,2}(\mathbb{R}^3).$$
We have to prove that
(3.2) $$(\Delta^{-1}u_k^2)u_k \longrightarrow hu \qquad \text{in } (H_r^1(\mathbb{R}^3))'.$$
Denote
$$q = \frac{12}{5}, \qquad r = \frac{12}{7}, \qquad \alpha = \frac{q}{r} = \frac{7}{5}, \qquad \alpha' = \frac{\alpha}{\alpha-1} = \frac{7}{2},$$
clearly
$$2 < q < 6, \qquad \frac{6}{5} < r < 2, \qquad \alpha' = \frac{6}{r}.$$

We have
$$\|(\Delta^{-1}u_k^2)u_k - hu\|_{L^r} \leq$$

(3.3) $$\leq \|(\Delta^{-1}u_k^2)u_k - (\Delta^{-1}u_k^2)u\|_{L^r} + \|(\Delta^{-1}u_k^2)u - hu\|_{L^r},$$

by Hölder inequality (note that $1/r = 1/6 + 1/q$)
$$\|(\Delta^{-1}u_k^2)u_k - (\Delta^{-1}u_k^2)u\|_{L^r} \leq$$
$$\leq \|\Delta^{-1}u_k^2\|_{L^6}\|u_k - u\|_{L^q},$$

then, using the compactness of the embedding (2.3), we see that $u_k \longrightarrow u$ in $L^q(\mathbb{R}^3)$ and $\{\Delta^{-1}u_k^2\}$ is bounded in $\mathcal{D}^{1,2}(\mathbb{R}^3)(\hookrightarrow L^6(\mathbb{R}^3))$, we have:

(3.4) $$\|(\Delta^{-1}u_k^2)u_k - (\Delta^{-1}u_k^2)u\|_{L^r} \longrightarrow 0.$$

From the fact that
$$u_k \longrightarrow u \text{ in } L^q(\mathbb{R}^3),$$
we see that
$$u_k^2 \longrightarrow u^2 \text{ in } L^{6/5}(\mathbb{R}^3).$$

Applying Lemma 2.1, we find
$$\Delta^{-1}u_k^2 \longrightarrow \Delta^{-1}u^2 \text{ in } \mathcal{D}^{1,2}(\mathbb{R}^3).$$

Now the Sobolev embedding (2.2) guarantees that
$$\Delta^{-1}u_k^2 \longrightarrow \Delta^{-1}u^2 \quad \text{in } L^6(\mathbb{R}^3).$$

Comparing this result with (3.1), we conclude that $h = \Delta^{-1}u^2$ and via
$$\|(\Delta^{-1}u_k^2)u - hu\|_{L^r} \leq \|\Delta^{-1}u_k^2 - h\|_{L^{\alpha'r}}\|u\|_{L^q},$$
we get

(3.5) $$\|(\Delta^{-1}u_k^2)u - hu\|_{L^r} \longrightarrow 0.$$

So we have, by (3.3), (3.4) and (3.5), that
$$(\Delta^{-1}u_k^2)u_k \longrightarrow hu \quad \text{in } L^r(\mathbb{R}^3).$$

From the properties (2.2) we arrive at (3.2).

We prove that $J_1$ is weakly continuous. Here it suffices to observe that the operator
$$Q: u \in H_r^1(\mathbb{R}^3) \longmapsto u^2 \in L^{6/5}(\mathbb{R}^3)$$
is compact, indeed, by the compact embeddings of $H_r^1(\mathbb{R}^3)$, the operator:
$$H_r^1(\mathbb{R}^3) \hookrightarrow\hookrightarrow L^{12/5}(\mathbb{R}^3) \xrightarrow{Q} L^{6/5}(\mathbb{R}^3)$$
is compact and, by Lemma 2.1, the following one
$$\Delta^{-1}: L^{6/5}(\mathbb{R}^3) \longrightarrow \mathcal{D}^{1,2}(\mathbb{R}^3)$$
is continuous.

We prove that $J_2$ is weakly continuous. Let $\{u_k\} \subset H^1(\mathbb{R}^3)$ and $u \in H^1(\mathbb{R}^3)$ such that
$$u_k \rightharpoonup u \quad \text{weakly in } H^1(\mathbb{R}^3).$$

Since
$$u_k \rightharpoonup u \quad \text{weakly in } L^2(\mathbb{R}^3),$$
there exists $C > 0$ such that
$$\|u_k\|_{L^2} \leq C, \qquad \|u\|_{L^2} \leq C, \qquad k \in \mathbb{N}.$$
Fix $\varepsilon > 0$ there results
$$(3.6) \quad \int_{\{|x|>z/\varepsilon\}} V(x)u_k^2 dx \leq C\varepsilon, \qquad \int_{\{|x|>z/\varepsilon\}} V(x)u^2 dx \leq C\varepsilon, \qquad k \in \mathbb{N}.$$
By the Sobolev inequality,
$$u_k^2 \rightharpoonup u^2 \quad \text{weakly in } L^3(\mathbb{R}^3),$$
since $V \in L^{\frac{3}{2}}(\{|x| \leq z/\varepsilon\})$, there results
$$\int_{\{|x|\leq z/\varepsilon\}} V(x)u_k^2 dx \longrightarrow \int_{\{|x|\leq z/\varepsilon\}} V(x)u^2 dx.$$
Then, by the previous one and (3.6), we can conclude
$$\int_{\mathbb{R}^3} V(x)u_k^2 dx \longrightarrow \int_{\mathbb{R}^3} V(x)u^2 dx.$$
Since, by well known arguments, the functional
$$u \in H^1(\mathbb{R}^3) \longmapsto \int_{\mathbb{R}^3} |\nabla u|^2 dx$$
is weakly lower semicontinuous, we are done. $\square$

**Lemma 3.2.** *The functional $J$ is coercive in $H_r^1(\mathbb{R}^3)$, i. e. for all sequence $\{u_k\} \subset H_r^1(\mathbb{R}^3)$ such that $\|u_k\|_{H^1} \longrightarrow +\infty$ there results $\lim\limits_{k} J(u_k) = +\infty$.*

PROOF. Denote
$$B_{H_r^1} = \{u \in H_r^1(\mathbb{R}^3) \big| \|u\|_{H^1} = 1\}.$$
Let $\{u_k\} \subset H_r^1(\mathbb{R}^3)$ be a sequence, such that
$$\|u_k\|_{H^1} \longrightarrow +\infty.$$
Take
$$\lambda_k = \|u_k\|_{H^1}$$
and set
$$\tilde{u}_k = \frac{u_k}{\lambda_k}.$$
Then obviously,
$$u_k = \lambda_k \tilde{u}_k$$
with $\lambda_k \in \mathbb{R}$ tending to $+\infty$ and $\tilde{u}_k \in B_{H_r^1}$. We have
$$J(u_k) = a_k \lambda_k^2 + b_k \lambda_k^4 - c_k \lambda_k^2,$$

with
$$a_k = \frac{1}{4}\int_{\mathbb{R}^3}|\nabla \tilde{u}_k|^2 dx \in \left[0, \frac{1}{4}\right],$$
$$b_k = \pi \int_{\mathbb{R}^3}|\nabla \Delta^{-1}\tilde{u}_k^2|^2 dx \geq 0,$$
$$c_k = \frac{1}{2}\int_{\mathbb{R}^3}V(x)\tilde{u}_k^2 dx \geq 0.$$

Observe that by Sobolev inequality there results

$$2c_k = \int_{\{|x|\leq 1\}}V(x)\tilde{u}_k^2 dx + \int_{\{|x|>1\}}V(x)\tilde{u}_k^2 dx \leq$$
$$\leq \|V\|_{L^{\frac{3}{2}}(\{|x|\leq 1\})}\|\tilde{u}_k\|_{L^6}^2 + \sup_{|x|\geq 1}V(x)\|\tilde{u}_k\|_{L^2}^2 \leq$$
$$\leq \left(C\|V\|_{L^{\frac{3}{2}}(\{|x|\leq 1\})} + \sup_{|x|\geq 1}V(x)\right)\|\tilde{u}_k\|_{H^1}^2 =$$
$$= \left(C\|V\|_{L^{\frac{3}{2}}(\{|x|\leq 1\})} + \sup_{|x|\geq 1}V(x)\right),$$

where $C > 0$ is the Sobolev embedding constant. Since, by Lemma 3.1, $u \in H^1_r(\mathbb{R}^3) \longmapsto \int_{\mathbb{R}^3}|\nabla \Delta^{-1}u^2|^2 dx$ is weakly continuous and $B_{H^1_r}$ is bounded in $H^1_r(\mathbb{R}^3)$, there exists $\alpha > 0$ such that $b_k \geq \alpha > 0$. Then we can conclude that

$$\lim_k J(u_k) = +\infty,$$

and so we are done. □

The two previous lemma guarantees that $J$ is bounded from below. Alternatively, we can give a direct proof of this fact.

**Lemma 3.3.** *The functional $J$ is bounded from below on $B$.*

PROOF. For each $u \in B$ there results

$$(3.7) \qquad J(u) \geq \frac{1}{4}\int_{\mathbb{R}^3}|\nabla u|^2 dx - \frac{1}{2}\int_{\mathbb{R}^3}V(x)|u|^2 dx.$$

By Kato's Inequality (see (7.13) page 35 in [15]) and since $u \in B$

$$\int_{\mathbb{R}^3}V(x)|u|^2 dx = z\int_{\mathbb{R}^3}\frac{|u|^2}{|x|}dx \leq Cz\|u\|_{L^2}\|\nabla u\|_{L^2} = CNz\|\nabla u\|_{L^2},$$

for some constant $C > 0$. So, by (3.7),

$$(3.8) \qquad J(u) \geq \frac{1}{4}\|\nabla u\|_{L^2}^2 - \frac{CNz}{2}\|\nabla u\|_{L^2}.$$

Since the map
$$x \in \mathbb{R} \longmapsto \frac{1}{4}x^2 - \frac{NCz}{2}x$$
is bounded from below, by (3.8), the claim is done. □

## 4. Spectral Results

The main result of this section is the following.

**Proposition 4.1.** *Let $(u, \omega) \in H^1_r(\mathbb{R}^3) \times \mathbb{R}$ be solution of the equation in (2.4). If*

(4.1) $$0 < \int_{\mathbb{R}^3} u^2 dx \leq N,$$

*then*

(4.2) $$\omega < 0.$$

The proposition implies that Theorem 1.3 is valid. To prove the above proposition some Lemmas are needed.

**Lemma 4.1.** *Let $u \in C^2(\{|x| \geq R\})$ be a solution of*

(4.3) $$\Delta u + p(x)u = 0, \qquad |x| \geq R,$$

*for some $R > 0$, if $p \in C(\mathbb{R}^3)$ and there exist $\alpha$, $R_0 > R$ such that*

(4.4) $$\frac{\partial p}{\partial r} + \frac{2(1-\alpha)}{|x|} p \geq 0, \qquad |x| \geq R_0,$$

*then*

(4.5) $$\liminf_{R \to +\infty} \frac{1}{R^\alpha} \int_{\{R_0 \leq |x| \leq R\}} p(x) u^2(x) dx > 0.$$

PROOF. It is direct consequence of [2, Theorem 3]. □

**Lemma 4.2.** *Let $u \in H^1_r(\mathbb{R}^3)$, $v \in L^1(\mathbb{R}^3) \cap L^{6/5}(\mathbb{R}^3)$ radial, $\omega \geq 0$ and*

(4.6) $$v \geq 0, \qquad \int_{\mathbb{R}^3} v dx \leq N$$

*and*

(4.7) $$-\frac{1}{2}\Delta u - 4\pi(\Delta^{-1}v)u - V(x)u = \omega u,$$

*then $u \equiv 0$.*

PROOF. Assume, by absurd, that there exist $u \not\equiv 0$ and $\omega \geq 0$ satisfying (4.6) and (4.7). Denote
$$p(x) \doteq 8\pi(\Delta^{-1}v)(x) + 2V(x) + 2\omega, \qquad x \in \mathbb{R}^3,$$
clearly $u$ is solution of (4.3). We shall apply Lemma 4.1 for this take $\alpha$, $0 < \alpha < \dfrac{1}{2}$. For $r \doteq |x|$, there results
$$\frac{\partial p}{\partial r}(x) + \frac{2(1-\alpha)}{|x|} p(x) =$$

(4.8) $$= 8\pi \left( \frac{\partial(\Delta^{-1}v)}{\partial r}(x) + \frac{2(1-\alpha)}{|x|}(\Delta^{-1}v)(x) \right) +$$
$$+ 2\left( \frac{\partial V}{\partial r}(x) + \frac{2(1-\alpha)}{|x|} V(x) \right) + \frac{4(1-\alpha)\omega}{|x|}.$$

Moreover, by (1.6),

$$(4.9) \qquad \frac{\partial V}{\partial r}(x) + \frac{2(1-\alpha)}{|x|}V(x) =$$

$$= -\frac{z}{r^2} + \frac{2(1-\alpha)z}{r^2} = \frac{(1-2\alpha)z}{r^2}.$$

By [16] or Lemma 8.2 in the Appendix,

$$(4.10) \qquad 4\pi(\Delta^{-1}v)(x) = -\int_{\mathbb{R}^3} \frac{v(y)}{\max\{|x|,|y|\}} dy, \qquad x \in \mathbb{R}^3,$$

so using this relation and Lemma 8.3 from the Appendix, we find

$$4\pi\left(\frac{\partial \Delta^{-1}v}{\partial r}(x) + \frac{2(1-\alpha)}{|x|}\Delta^{-1}v(x)\right) =$$

$$= \int_{|y|<r} \frac{v(y)}{|x|^2} dy - \frac{2(1-\alpha)}{r}\int_{\mathbb{R}^3} \frac{v(y)}{\max\{|x|,|y|\}} dy =$$

$$= \int_{|y|<r} \frac{v(y)}{|x|^2} dy - \frac{2(1-\alpha)}{r}\int_{\{|y|\leq r\}} \frac{v(y)}{\max\{|x|,|y|\}} dy -$$

$$- \frac{2(1-\alpha)}{r}\int_{\{|y|\geq r\}} \frac{v(y)}{\max\{|x|,|y|\}} dy =$$

$$= \int_{\{|y|\leq r\}} \frac{v(y)}{\max\{|x|^2,|y|^2\}} dy - \frac{2(1-\alpha)}{r^2}\int_{\{|y|\leq r\}} v(y)dy -$$

$$(4.11) \quad - \frac{2(1-\alpha)}{r}\int_{\{|y|\geq r\}} \frac{v(y)}{|y|} dy \geq$$

$$\geq \int_{\{|y|\leq r\}} \frac{v(y)}{r^2} dy - \frac{2(1-\alpha)}{r^2}\int_{\{|y|\leq r\}} v(y)dy - \frac{2(1-\alpha)}{r^2}\int_{\{|y|\geq r\}} v(y)dy \geq$$

$$\geq -\frac{(1-2\alpha)}{r^2}\int_{\mathbb{R}^3} v(y)dy - \frac{2(1-\alpha)}{r^2}\int_{\{|y|\geq r\}} v(y)dy.$$

By (4.8), (4.9) and (4.11),

$$\frac{\partial p}{\partial r}(x) + \frac{2(1-\alpha)}{|x|}p(x) \geq$$

$$(4.12) \qquad \geq 2\frac{(1-2\alpha)}{r^2}\left(z - \int_{\mathbb{R}^3} v(y)dy\right) +$$

$$+ 4\frac{(1-\alpha)}{r}\left(\omega - \frac{1}{r}\int_{\{|y|\geq r\}} v(y)dy\right).$$

There exists $R_0 > 0$ such that
$$\frac{1}{|x|} \int_{\{|y| \geq |x|\}} v(y) dy \leq \omega, \qquad |x| \geq R_0,$$

by (1.8), (4.6) and (4.12), since $0 < \alpha < \frac{1}{2}$, we have

(4.13) $$\frac{\partial p}{\partial r}(x) + \frac{2(1-\alpha)}{|x|} p(x) \geq 0, \qquad |x| \geq R_0.$$

By (4.3) and Lemma 4.1, the formula (4.5) holds true.

On the other hand, we have

(4.14) $$\int_{\mathbb{R}^3} u^2 (\Delta^{-1} v) dx \leq \|u\|_{L^{12/5}}^2 \|\Delta^{-1} v\|_{L^6}$$

and, as in Lemma 3.2,

(4.15) $$\int_{\mathbb{R}^3} u^2 V dx \leq \|V\|_{L^{3/2}(\{|x| \leq 1\})} \|u\|_{L^6}^2 + z \|u\|_{L^2}^2,$$

so, by (4.14) and (4.15),

$$\int_{\{R_0 \leq |x| \leq R\}} pu^2 dx \leq \int_{\mathbb{R}^3} pu^2 dx =$$

(4.16) $$= 2 \left( 4\pi \int_{\mathbb{R}^3} (\Delta^{-1} v) u^2 dx + \int_{\mathbb{R}^3} V u^2 dx + \int_{\mathbb{R}^3} \omega u^2 dx \right) \leq$$

$$\leq 8\pi \|u\|_{L^{12/5}}^2 \|\Delta^{-1} v\|_{L^6} + 2\|V\|_{L^{3/2}(\{|x| \leq 1\})} \|u\|_{L^6}^2 + 2z \|u\|_{L^2}^2 + 2\omega \|u\|_{L^2}^2.$$

Then

(4.17) $$\lim_{R \to +\infty} \frac{1}{R^\alpha} \int_{\{R_0 \leq |x| \leq R\}} p(x) u^2(x) dx = 0,$$

and this is absurd, since (4.17) contradicts (4.5), this concludes the proof. □

**Corollary 4.1.** *If $V \equiv 0$ and the assumptions of Lemma 4.2 are satisfied, then $u \equiv 0$.*

PROOF. Suppose, by absurd, that there is $u \not\equiv 0$ solution of (4.7), multiplying by $u$ and integrating on $\mathbb{R}^3$, we get $\omega > 0$. We are going to apply the Agmon's result of Lemma 4.1. For this we have to verify the condition (4.4) for $|x|$ large enough, $0 < \alpha < \frac{1}{2}$ and
$$p(x) \doteq 2(\Delta^{-1} v)(x) + 2\omega, \qquad x \in \mathbb{R}^3.$$

The argument of the previous lemma (with $z = 0$) gives

$$\frac{\partial p}{\partial r}(x) + \frac{2(1-\alpha)}{|x|} p(x) \geq$$

$$\geq 4 \frac{(1-\alpha)\omega}{r} - 2 \frac{(1-2\alpha)}{r^2} \int_{\mathbb{R}^3} v(y) dy - 4 \frac{(1-\alpha)}{r^2} \int_{\{|y| \geq r\}} v(y) dy.$$

So, for $R_0 > 0$ sufficiently large

$$\frac{\partial p}{\partial r}(x) + \frac{2(1-\alpha)}{|x|} p(x) \geq 0, \qquad |x| \geq R_0.$$

By (4.16), with $z = 0$ we have

$$\frac{1}{R^\alpha} \int_{\{R_0 \leq |x| \leq R\}} pu^2 dx \leq \frac{2}{R^\alpha} \left( \|u\|_{L^{12/5}}^2 \|\Delta^{-1} v\|_{L^6} + \omega \|u\|_{L^2}^2 \right).$$

This is absurd, because it contradicts (4.5), then $u \equiv 0$. □

PROOF OF THEOREM 1.1. Denote

$$v(x) \doteq u^2(x), \qquad x \in \mathbb{R}^3.$$

By the Sobolev inequalities

$$v \in L^1(\mathbb{R}^3) \cap L^{6/5}(\mathbb{R}^3),$$

it is radial and, by the constraint in (2.4), it satisfies also (4.6). Since $\omega \geq 0$, the claim is direct consequence of the previous corollary and of the equivalence between (1.7) and (2.4). □

**Lemma 4.3.** *Let $u \in H_r^1(\mathbb{R}^3)$ be a solution of the equation in (2.4) such that*

(4.18) $$\int_{\mathbb{R}^3} u^2 dx \leq N,$$

*if*

(4.19) $$\omega \geq 0,$$

*then $u \equiv 0$.*

PROOF. Denote

$$v(x) \doteq u^2(x), \qquad x \in \mathbb{R}^3.$$

By the Sobolev inequalities

$$v \in L^1(\mathbb{R}^3) \cap L^r(\mathbb{R}^3), \qquad \frac{6}{5} < r \leq 2,$$

it is radial and, by (4.18), it satisfies also (4.6). Since $\omega \geq 0$, the claim is direct consequence of the previous lemma. □

PROOF OF PROPOSITION 4.1. It is direct consequence of the previous lemma. □

PROOF OF THEOREM 1.3. It is direct consequence of Proposition 4.1 and of the equivalence between (1.7) and (2.4). □

## 5. PROOF OF THEOREM 1.2

In this section we shall prove Theorem 1.2 We begin proving some lemmas.

**Lemma 5.1.** *The functional $J\big|_{H^1_r(\mathbb{R}^3)\cap B}$ satisfies the Palais-Smale condition in each level $]-\infty, -\varepsilon]$, $\varepsilon > 0$, i.e. any sequence $\{u_k\} \subset H^1_r(\mathbb{R}^3) \cap B$ such that $\{J(u_k)\}$ is bounded and*

$$\tag{5.1} J(u_k) \leq -\varepsilon, \qquad J'\big|_{H^1_r(\mathbb{R}^3)\cap B}(u_k) \longrightarrow 0,$$

*contains a converging subsequence.*

PROOF. Fix $\varepsilon > 0$. Let $\{u_k\} \subset H^1_r(\mathbb{R}^3) \cap B$ be such that $\{J(u_k)\}$ is bounded and satisfies (5.1). First of all observe that, by iii) of Lemma 2.2, there results

$$J'\big|_{H^1_r(\mathbb{R}^3)\cap B}(u) = 0 \iff J'\big|_B(u) = 0,$$

then we can suppose

$$J'\big|_B(u_k) \longrightarrow 0.$$

Since $J(u_k) \leq -\varepsilon$, by Lemma 3.2, $\{u_k\}$ is bounded in $H^1_r(\mathbb{R}^3)$, passing to a subsequence, there exists $u \in H^1_r(\mathbb{R}^3)$ such that

$$\tag{5.2} u_k \rightharpoonup u \qquad \text{weakly in } H^1_r(\mathbb{R}^3).$$

We shall prove that

$$\tag{5.3} u_k \longrightarrow u \qquad \text{in } H^1_r(\mathbb{R}^3).$$

By definition, there exists $\{\omega_k\} \subset \mathbb{R}$ such that

$$J'\big|_B(u_k) = J'(u_k) - \omega_k u_k, \qquad k \in \mathbb{N}.$$

Observe that, since $\{u_k\} \subset B$, we have

$$N\omega_k = \langle J'\big|_B(u_k), u_k \rangle - \langle J'(u_k), u_k \rangle =$$
$$= \langle J'\big|_B(u_k), u_k \rangle - \frac{1}{2}\int_{\mathbb{R}^3} |\nabla u_k|^2 dx -$$
$$- 4\pi \int_{\mathbb{R}^3} |\nabla \Delta^{-1} u_k^2|^2 dx + \int_{\mathbb{R}^3} V(x)|u_k|^2 dx =$$
$$= \langle J'\big|_B(u_k), u_k \rangle - 2J(u_k) - 2\pi \int_{\mathbb{R}^3} |\nabla \Delta^{-1} u_k^2|^2 dx,$$

by Lemma 3.1 and (5.1), $\{\omega_k\}$ is bounded in $\mathbb{R}$ and so passing to a subsequence there results

$$\tag{5.4} \omega_k \longrightarrow \omega$$

and

$$\tag{5.5} -\frac{1}{2}\Delta u - 4\pi(\Delta^{-1} u^2)u - V(x)u = \omega u.$$

If $\omega < 0$, by Lemma 3.1, (5.2), (5.4) and (5.5),

$$\frac{1}{2}\int_{\mathbb{R}^3}|\nabla u_k|^2 dx - \omega \int_{\mathbb{R}^3} u_k^2 dx =$$
$$= \langle J\,|'_B(u_k), u_k\rangle - 4\pi \int_{\mathbb{R}^3}|\nabla\Delta^{-1}u_k^2|^2 dx +$$
$$+ \int_{\mathbb{R}^3} V(x)u_k^2 dx + (\omega_k - \omega)\int_{\mathbb{R}^3} u_k^2 dx \longrightarrow$$
$$\longrightarrow -4\pi \int_{\mathbb{R}^3}|\nabla\Delta^{-1}u^2|^2 dx + \int_{\mathbb{R}^3} V(x)u^2 dx =$$
$$= \frac{1}{2}\int_{\mathbb{R}^3}|\nabla u|^2 dx - \omega \int_{\mathbb{R}^3} u^2 dx,$$

and then (5.3).

Now we consider the case $\omega \geq 0$. If $\|u\|_{L^2} = 0$, by Lemma 3.1, we have

$$0 = J(u) \leq \liminf_k J(u_k) \leq -\varepsilon,$$

that is absurd. If $0 < \|u\|_{L^2} \leq N$ then $u$ is solution of the equation in (2.4), (4.18) and (4.19) hold. So, by Lemma 4.3, we have $u \equiv 0$ and also this is absurd. This concludes the proof. $\square$

*Remark* 5.1. Let $\rho \in L^1(\mathbb{R}^3) \cap L^r(\mathbb{R}^3)$, with $\frac{6}{5} < r \leq 2$, $\varphi \in \mathcal{D}^{1,2}(\mathbb{R}^3)$ radially symmetric maps such that

$$\Delta\varphi = \rho.$$

Denote

$$\rho_\nu(x) \doteq \rho(\nu x), \qquad x \in \mathbb{R}^3, \ \nu \geq 0$$

we claim that the unique solution $\varphi_\nu$ of the equation

$$\Delta\varphi_\nu = \rho_\nu$$

is

$$\varphi_\nu(x) \doteq \nu^{-2}\varphi(\nu x), \qquad x \in \mathbb{R}^3.$$

Indeed, denoting $r \doteq |x|$, there results

$$\Delta\varphi_\nu(x) = \Delta\varphi_\nu(r) = \partial^2_{rr}\varphi_\nu(r) + \frac{2}{r}\partial_r\varphi_\nu(r) =$$
$$= \frac{1}{\nu^2}\left(\nu^2\partial^2_{rr}\varphi(\nu r) + \frac{2\nu}{r}\partial_r\varphi(\nu r)\right) =$$
$$= \partial^2_{rr}\varphi(\nu r) + \frac{2}{r\nu}\partial_r\varphi(\nu r) =$$
(5.6) $$= \Delta\varphi(\nu x) = \rho(\nu x) = \rho_\nu(x).$$

Moreover

$$\int_{\mathbb{R}^3}|\nabla\Delta^{-1}\rho_\nu(x)|^2 dx = \int_{\mathbb{R}^3}|\nabla\varphi_\nu(x)|^2 dx =$$

(5.7) $$= \frac{1}{\nu^2}\int_{\mathbb{R}^3}|\nabla\varphi(\nu x)|^2 dx = \frac{1}{\nu^5}\int_{\mathbb{R}^3}|\nabla\varphi(x)|^2 dx = \frac{1}{\nu^5}\int_{\mathbb{R}^3}|\nabla\Delta^{-1}\rho(x)|^2 dx.$$

In the last part of this section we need some more notations. Define

$$c_k \doteq \inf \{ \sup J(A) | A \in \mathcal{A},\ \gamma(A) \geq k \}, \qquad k \in \mathbb{N}\setminus\{0\},$$
$$\tilde{c}_k \doteq \inf \{ \sup J(h(S^{k-1})) | h \in \Omega_k \}, \qquad k \in \mathbb{N}\setminus\{0\},$$
$$\tilde{c}_{k,\lambda} \doteq \inf \{ \sup J(h(S^{k-1})) | h \in \Omega_{k,\lambda} \}, \qquad k \in \mathbb{N}\setminus\{0\},\ \lambda > 0,$$

where

$$\mathcal{A} \doteq \{ A \subset H^1_r(\mathbb{R}^3) \cap B | A \text{ closed and symmetric} \},$$
$$\Omega_k \doteq \{ h : S^{k-1} \longrightarrow H^1_r(\mathbb{R}^3) \cap B | h \text{ continuous and odd} \}, \qquad k \in \mathbb{N}\setminus\{0\},$$
$$\Omega_{k,\lambda} \doteq \{ h : S^{k-1} \longrightarrow H^1_r(\mathbb{R}^3) \cap B_\lambda | h \text{ continuous and odd} \}, \qquad k \in \mathbb{N}\setminus\{0\},\ \lambda > 0,$$
$$B_\lambda \doteq \{ u \in H^1(\mathbb{R}^3) | \|u\|_{L^2} = \lambda \}, \qquad \lambda > 0$$

and $\gamma$ is the Genus (see e. g. [3, Definition 1.1]).

**Lemma 5.2.** *There results*

(5.8) $$c_k \leq \tilde{c}_k \leq \tilde{c}_{k,\lambda},$$

*for each $k \in \mathbb{N}\setminus\{0\}$ and $0 < \lambda \leq \sqrt{N}$.*

PROOF. Fix $k \in \mathbb{N}\setminus\{0\}$. We prove that

(5.9) $$c_k \leq \tilde{c}_k.$$

Let $h \in \Omega_k$, since $h$ is continuous and odd the set $J(h(S^{k-1}))$ is closed and symmetric. Moreover $h(S^{k-1}) \subset B$ and, by the invariance property of the Genus, there results

$$\gamma(h(S^{k-1})) \geq \gamma(S^{k-1}) = k.$$

So we have

$$c_k \leq \sup J(h(S^{k-1}))$$

and then (5.9) is proved.

We prove that

(5.10) $$\tilde{c}_k \leq \tilde{c}_{k,\lambda} \qquad 0 < \lambda \leq \sqrt{N}.$$

Fix $0 < \lambda \leq \sqrt{N}$ and define

$$h_\lambda(\xi)(x) = \frac{1}{\lambda^5} h(\xi)\Big(\frac{x}{\lambda^4}\Big), \qquad h \in \Omega_k,\ \xi \in S^{k-1}.$$

Let $h \in \Omega_k$ and $\xi \in S^{k-1}$ such that

(5.11) $$\frac{3}{2N^2} \int_{\mathbb{R}^3} |\nabla u|^2 dx - \int_{\mathbb{R}^3} V(x)|u|^2 dx \geq 0,$$

where $u \doteq h(\xi)$. Set

$$\nu \doteq \frac{1}{\lambda^4}, \qquad u_\nu(x) \doteq h_\lambda(x) = \nu^{\frac{5}{4}}(x) u(\nu x)$$

and observe that, by (5.7), there results

$$\int_{\mathbb{R}^3} |u_\nu|^2 dx = \frac{1}{\nu^{1/2}} \int_{\mathbb{R}^3} |u|^2 dx = \lambda^2 N,$$

$$\int_{\mathbb{R}^3} |\nabla u_\nu|^2 dx = \nu^{\frac{3}{2}} \int_{\mathbb{R}^3} |\nabla u|^2 dx,$$

$$\int_{\mathbb{R}^3} |\nabla \Delta^{-1} u_\nu^2(x)|^2 dx = \int_{\mathbb{R}^3} |\nabla \Delta^{-1} u^2(x)|^2 dx,$$

$$\int_{\mathbb{R}^3} V(x) |u_\nu|^2 dx = \nu^{\frac{1}{2}} \int_{\mathbb{R}^3} V(x) |u|^2 dx.$$

Consider the map

$$f(\nu) \doteq J(u_\nu) = \frac{\nu^{\frac{3}{2}}}{4} \int_{\mathbb{R}^3} |\nabla u|^2 dx + \pi \int_{\mathbb{R}^3} |\nabla \Delta^{-1} u^2|^2 dx - \frac{\nu^{\frac{1}{2}}}{2} \int_{\mathbb{R}^3} V(x) |u|^2 dx,$$

there results

$$\frac{df}{d\nu}(\nu) = \frac{3\nu^{\frac{1}{2}}}{8} \int_{\mathbb{R}^3} |\nabla u|^2 dx - \frac{1}{4\nu^{\frac{1}{2}}} \int_{\mathbb{R}^3} V(x) |u|^2 dx.$$

Clearly

$$\frac{df}{d\nu}(\nu) \geq 0 \iff \frac{3\nu}{2} \int_{\mathbb{R}^3} |\nabla u|^2 dx - \int_{\mathbb{R}^3} V(x) |u|^2 dx \geq 0$$

and then, by (5.11), $f$ is increasing for $\nu \geq 1/N^2$, namely

$$J(h(\xi)) = J(u) \leq J(u_\nu) = J(h_\lambda(\xi)).$$

Since, if there exists $\xi' \in S^{k-1}$, $\xi \neq \xi'$ such that $h(\xi')$ does not satisfy (5.11), we have $J(h(\xi')) \leq J(h(\xi))$, then

$$\sup J\big(h\big(S^{k-1}\big)\big) \leq \sup J\big(h_\lambda\big(S^{k-1}\big)\big).$$

This concludes the proof of (5.10). □

**Lemma 5.3.** *For all $k \in \mathbb{N}\setminus\{0\}$, there exist a subspace $V_k \subset H^1_r(\mathbb{R}^3)$ of dimension $k$ and $\nu > 0$ such that*

$$\int_{\mathbb{R}^3} \Big(\frac{1}{2}|\nabla u|^2 - V(x) u^2\Big) dx \leq -\nu,$$

*for all $u \in V_k \cap B$.*

PROOF. Let $u$ be a smooth map with compact support such that

$$\int_{\mathbb{R}^3} |u|^2 dx = N, \quad \text{supp}(u) \subset B_2(0) \setminus B_1(0),$$

where

$$B_\rho(x) \doteq \big\{ y \in \mathbb{R}^3 \,\big|\, |x - y| < \rho \big\}, \qquad x \in \mathbb{R}^3, \ \rho > 0.$$

Denote

$$u_\lambda(x) \doteq \lambda^{\frac{3}{2}} u(\lambda x), \qquad \lambda > 0, \ x \in \mathbb{R}^3,$$

there results

$$\int_{\mathbb{R}^3} |u|^2 dx = \int_{\mathbb{R}^3} |u_\lambda|^2 dx = N, \quad \text{supp}(u_\lambda) \subset B_{\frac{2}{\lambda}}(0) \setminus B_{\frac{1}{\lambda}}(0).$$

We have
$$\int_{\mathbb{R}^3} \left(\frac{1}{2}|\nabla u_\lambda|^2 - V(x)u_\lambda^2\right)dx = \int_{\mathbb{R}^3} \left(\lambda^2 \frac{1}{2}|\nabla u|^2 - V(\frac{x}{\lambda})u^2\right)dx \le$$
$$\le \lambda^2 \int_{\mathbb{R}^3} \frac{1}{2}|\nabla u|^2 dx - N \inf_{\frac{\operatorname{supp} u}{\lambda}} V \le \lambda^2 \int_{\mathbb{R}^3} \frac{1}{2}|\nabla u|^2 dx - \frac{z\lambda}{2}N.$$

There exists $\lambda_0 > 0$ such that
$$\int_{\mathbb{R}^3} \left(\frac{1}{2}|\nabla u_{\lambda_0}|^2 - V(x)u_{\lambda_0}^2\right)dx < 0.$$

Let $k \in \mathbb{N}\backslash\{0\}$ and $u_1, u_2, \ldots, u_k$ smooth maps with compact supports such that
$$\int_{\mathbb{R}^3} |u_i|^2 dx = 1, \quad \operatorname{supp}(u_i) \subset B_{2i}(0)\backslash B_i(0), \quad i = 1, 2, \ldots, k.$$

Using an analogous argument we are able to find $\lambda_1, \lambda_2, \ldots, \lambda_k > 0$ such that
$$\int_{\mathbb{R}^3} \left(\frac{1}{2}|\nabla u_{i_{\lambda_i}}|^2 - V(x)u_{i_{\lambda_i}}^2\right)dx < 0, \quad i = 1, 2, \ldots, k.$$

Let
$$0 < \bar{\lambda} < \min\{\lambda_1, \lambda_2, \ldots, \lambda_k\}$$
and $V_k$ the subspace spanned by $u_{1_{\bar\lambda}}, u_{2_{\bar\lambda}}, \ldots, u_{k_{\bar\lambda}}$. Since the supports of this maps are pairwise disjoint, $V_k$ has dimension $k$. Since for all $i = 1, 2, \ldots, k$ and $\lambda \le \lambda_i$, there results
$$\int_{\mathbb{R}^3} \left(\frac{1}{2}|\nabla u_{i_\lambda}|^2 - V(x)u_{i_\lambda}^2\right) < 0$$
and $V_k \cap B$ is compact, the claim is proved. $\square$

**Lemma 5.4.** *There results*

(5.12) $$c_k < 0,$$

*for each $k \in \mathbb{N}\backslash\{0\}$.*

PROOF. Let $k \in \mathbb{N}\backslash\{0\}$, by Lemma 5.3 there exist $V_k \subset H^1_r(\mathbb{R}^3)$ subspace of dimension $k$ and $\nu > 0$ such that, for all $u \in V_k \cap B$,
$$\int_{\mathbb{R}^3} \left(\frac{1}{2}|\nabla u|^2 - V(x)u^2\right)dx \le -\nu.$$

Let $\lambda > 0$ and define
$$h_\lambda : V_k \cap B \longrightarrow H^1_r(\mathbb{R}^3), \quad h_\lambda(u) = \lambda^{\frac{1}{2}}u.$$

Fixed $u \in V_k \cap B$ and $0 < \lambda < \sqrt{N}$, there results

(5.13) $$J(h_\lambda(u)) \le -\frac{\lambda}{2}\nu + c\lambda^2 \le -\frac{\lambda}{2}\nu + c\lambda^2,$$

where $c$ is a positive constant. Then there exists $0 < \bar\lambda < \sqrt{N}$ such that for all $u \in V_k \cap B$ there results $J(h_{\bar\lambda}(u)) < 0$. Since $h_{\bar\lambda} \in \Omega_{\bar\lambda}$ and $V_k \cap B \simeq S^{k-1}$, by Lemma 5.2 and the compactness of $S^{k-1}$, we have
$$c_k \le \tilde c_k \le \tilde c_{k,\bar\lambda} \le \sup J(h_{\bar\lambda}(V_k \cap B)) < 0,$$

and we are done. $\square$

**Corollary 5.1.** *There results*

(5.14) $$\inf_{u \in H^1_r(\mathbb{R}^3) \cap B} J(u) < 0.$$

PROOF. It is direct consequence of the previous Lemma. □

**Lemma 5.5.** *Let $k \in \mathbb{N}$, $E \subset H^1(\mathbb{R}^3)$ be a subspace of dimension $k$ and $A \in \mathcal{A}$, if*

(5.15) $$\gamma(A) \geq k+1$$

*then*

(5.16) $$A \cap E^\perp \neq \emptyset.$$

PROOF. Assume, by absurd that (5.16) is false, there results

(5.17) $$P(A) \subset E \backslash \{0\},$$

where $P : H^1(\mathbb{R}^3) \longrightarrow E$ is the orthogonal projection on $E$. So we have

(5.18) $$\gamma(P(A)) \leq k.$$

On the other side, since $P$ is continuous and odd, by the invariance property of the Genus there results

$$k + 1 \leq \gamma(A) \leq \gamma(P(A)).$$

Since this is in contradiction with (5.18) the proof is done. □

**Lemma 5.6.** *The functional $J$ has a sequence $\{u_k\}_{k \in \mathbb{N}} \subset H^1_r(\mathbb{R}^3) \cap B$ of critical points such that*

$$\omega_k < 0, \qquad \omega_k \longrightarrow 0,$$

*where $\{\omega_k\}_{k \in \mathbb{N}} \subset \mathbb{R}$ is the sequence of the Lagrange multipliers associated to the critical points.*

PROOF. By Lemmas 5.1 and 5.4 (see [18, Theorem 9.1]) there exists a sequence $\{u_k\}_{k \in \mathbb{N}} \subset H^1_r(\mathbb{R}^3) \cap B$ of critical points of the functional $J$. Call $\{\omega_k\}_{k \in \mathbb{N}} \subset \mathbb{R}$ the sequence of the Lagrange multipliers associated to this critical points, namely

$$J'(u_k) - \omega_k u_k = 0, \qquad k \in \mathbb{N} \backslash \{0\}.$$

By Lemma 4.2, there results

$$\omega_k < 0, \qquad k \in \mathbb{N} \backslash \{0\}.$$

We have to prove that

(5.19) $$\omega_k \longrightarrow 0.$$

Let $\{V_k\}$ be a sequence of subspaces of $H^1_r(\mathbb{R}^3)$, such that

$$\dim(V_k) = k, \qquad \bigcup_{k \in \mathbb{N} \backslash \{0\}} V_k \text{ is dense in } H^1_r(\mathbb{R}^3).$$

Moreover, let $\{A_k\} \subset \mathcal{A}$ such that

(5.20) $$\gamma(A_k) \geq k, \qquad c_k \leq \sup J(A_k) \leq \frac{c_k}{2}, \qquad k \in \mathbb{N} \backslash \{0\}.$$

Call

$$W_k \doteq V^\perp_{k-1}, \qquad k \in \mathbb{N} \backslash \{0\},$$

unused20       GIUSEPPE MARIA COCLITE AND VLADIMIR GEORGIEV

by Lemma 5.5, there results

$$W_k \cap A_k \neq \emptyset, \quad k \in \mathbb{N}\setminus\{0\}.$$

Let $\{v_k\} \subset H^1_r(\mathbb{R}^3) \cap B$ such that

$$v_k \in W_k \cap A_k, \quad k \in \mathbb{N}\setminus\{0\},$$

clearly

(5.21) $\qquad v_k \rightharpoonup 0 \quad \text{weakly in } H^1_r(\mathbb{R}^3)$

and, by (5.20),

(5.22) $\qquad \sup J(V_k) \leq \dfrac{c_k}{2}, \quad k \in \mathbb{N}\setminus\{0\}.$

By (5.21) and Lemma 3.1 we have

(5.23) $\qquad 0 \leq \liminf_k J(v_k)$

and, by (5.22),

(5.24) $\qquad \limsup_k J(v_k) \leq \lim_k \dfrac{c_k}{2} \leq 0.$

By (5.23) and (5.24), we deduce $c_k \longrightarrow 0$. Since $2c_k \leq \omega_k < 0$, (5.19) is done. $\square$

PROOF OF THEOREM 1.2. Since

$$F(u, 4\pi\Delta^{-1}u^2) = J(u)$$

for all $u \in H^1(\mathbb{R}^3)$, by Lemma 2.2 and the previous one the claim is done. $\square$

## 6. PROOF OF THEOREM 1.4

Our next step is to show that the radially symmetric solutions

(6.1) $\qquad u \in H^1(\mathbb{R}^3), \qquad \nabla\varphi \in L^2(\mathbb{R}^3),$

to

$$-\frac{1}{2}\Delta u - \varphi u - \frac{z}{|x|}u = \omega u, x \in \mathbb{R}^3,$$

$$\Delta\varphi = 4\pi u^2, \; x \in \mathbb{R}^3,$$

(6.2) $\qquad \displaystyle\int_{\mathbb{R}^3} u^2 dx = N,$

constructed in the previous section, are more regular. More precisely, we shall derive the higher regularity

(6.3) $\qquad \nabla u \in H^{k-1}(|x| > \varepsilon), \nabla\varphi \in H^{k-1}(|x| > \varepsilon),$

where $k$ is arbitrary integer and $\varepsilon > 0$.

**Lemma 6.1.** *If the assumption* (6.1) *is satisfied, then*

(6.4) $\qquad \nabla u \in H^1(\mathbb{R}^3), \nabla\varphi \in H^1(\mathbb{R}^3), u \in L^\infty, \varphi \in L^\infty.$

PROOF. The assumption (6.1) and the Sobolev embedding in $\mathbb{R}^3$ guarantee that

(6.5) $$\varphi \in L^6, \ u \in L^p, 2 \leq p \leq 6.$$

This property and the Hölder inequality imply that the nonlinear term $\varphi u$ in the first equation in (6.2) is in $L^2$. The fact that $|x|^{-1}u \in L^2$ follows from the Hardy inequality and the fact that $\nabla u \in L^2$. Therefore this equation shows that $\Delta u \in L^2$, so $u \in H^2$. Using the second equation and the fact that $u^2 \in L^2$ we conclude that $\nabla \varphi \in H^1$. Finally the property $u \in L^\infty$ follows from the estimate

$$\|u\|_{L^\infty} \leq C \|\nabla u\|_{H^1(\mathbb{R}^3)}.$$

This estimate follows from the Fourier representation

$$u(x) = (2\pi)^{-3} \int_{\mathbb{R}^3} e^{-ix\xi} \hat{u}(\xi) d\xi,$$

the Cauchy inequality and the fact that

$$|\xi|^{-1}(1+|\xi|)^{-1} \in L^2(\mathbb{R}^3).$$

The Lemma is established. □

In the same way, proceeding inductively, we obtain the following.

**Lemma 6.2.** *If the assumption* (6.1) *is satisfied, then for any integer $k \geq 2$ and for any positive number $\varepsilon > 0$ we have*

(6.6) $$u \in H^k(|x| > \varepsilon), \qquad \nabla \varphi \in H^{k-1}(|x| > \varepsilon).$$

To study more precisely the behavior of the solution $u(x) = u(|x|)$ we introduce polar coordinates $r = |x|$ and set

(6.7) $$\mathbb{U}(r) = ru(r), \qquad \mathbb{V}(r) = -r\varphi(r).$$

Using the identities

$$\Delta \left( \frac{\mathbb{U}(r)}{r} \right) = \frac{\mathbb{U}''(r)}{r},$$

$$\Delta \left( \frac{\mathbb{V}(r)}{r} \right) = \frac{\mathbb{V}''(r)}{r},$$

where $\mathbb{U}'(r) = \partial_r \mathbb{U}(r)$, we can rewrite (6.2) in the form

$$-\frac{\mathbb{U}''}{2} + \frac{\mathbb{V}}{r}\mathbb{U} - \frac{z}{r}\mathbb{U} = \omega \mathbb{U}, r > 0,$$

(6.8) $$-\mathbb{V}'' = 4\pi \frac{\mathbb{U}^2}{r}, \ r > 0.$$

We shall need the following

**Lemma 6.3.** *Let $k \geq 1$ be an integer and $\varepsilon > 0$ be a real number. We have the following properties:*
  *a) if $u(x) = u(|x|) \in H^k(\mathbb{R}^3)$, then $\mathbb{U}(r) \in H^k(0, \infty)$;*
  *b) $u(x) = u(|x|) \in H^k(|x| > \varepsilon)$, if and only if $\mathbb{U}(r) \in H^k(\varepsilon, \infty)$.*

PROOF. The proof of a) follows from the relation

$$\partial_r^k \mathbb{U}(r) = r \partial_r^k u(r) + k \partial_r^{k-1} u(r)$$

valid for any integer $k \geq 1$. Note that the Hardy inequality implies
$$\int_0^\infty |\partial_r^{k-1} u(r)|^2 dr \leq C \|u\|_{H^k(\mathbb{R}^3)}^2.$$

For the property b) we can use the relation
$$\partial_r^k u(r) = \sum_{j=1}^k \frac{c_{k,j}}{r^j} \partial_r^{k-j} \mathbb{U}(r)$$

and the fact that $r^{-j}$ is bounded for $r \geq \varepsilon > 0$. $\square$

**Lemma 6.4.** *The functions $\mathbb{U}(r), \mathbb{V}(r)$ are smooth near $r = 0$.*

PROOF. From $u \in H^2$ (see Lemma 6.1) it follows $u \in L^\infty$, so
$$|\mathbb{U}(r)| = r|u(r)| \leq Cr$$
near $r = 0$. In the same way $\varphi \in L^\infty$ (Lemma 6.1) implies that
$$|\mathbb{V}(r)| = r|\varphi(r)| \leq Cr$$
near $r = 0$.

The system (6.2) shows that
$$|\mathbb{U}''(r)| + |\mathbb{V}''(r)| \leq C$$
so $\mathbb{U}(r), \mathbb{V}(r) \in C^1([0,1])$. Setting
$$a_1 = \mathbb{U}'(0), b_1 = \mathbb{V}'(0),$$
we can make the representation
$$\mathbb{U}(r) = a_1 r + \mathbb{U}_1(r), \ \mathbb{V}(r) = b_1 r + \mathbb{V}_1(r),$$
where $\mathbb{U}_1, \mathbb{V}_1 \in o(r)$ satisfy
$$-\frac{\mathbb{U}_1''}{2} + \frac{\mathbb{V}_1}{r} \mathbb{U}_1 - \frac{z}{r} \mathbb{U}_1 - \omega \mathbb{U}_1 = c_1 + O(r), r > 0,$$
$$(6.9) \qquad -\mathbb{V}_1'' - 4\pi \frac{\mathbb{U}_1^2}{r} = O(r), \ r > 0,$$
where $c_1 = \omega a_1$. These equations imply
$$\mathbb{U}_1''(r) = c_1 + O(r), \mathbb{V}_1''(r) = O(r),$$
so
$$\mathbb{U}_1(r) = \frac{c_1 r^2}{2} + O(r^3), \mathbb{V}_1(r) = O(r^3)$$
near $r = 0$ and these relations imply
$$\mathbb{U}_1(r), \mathbb{V}_1(r) \in C^2([0,1]).$$

Continuing further we obtain inductively
$$\mathbb{U}(r) = a_1 r + a_2 r^2 + \cdots + a_k r^k + \mathbb{U}_k(r), \ \mathbb{V}(r) = b_1 r + b_2 r^2 + \cdots + b_k r^k + \mathbb{V}_k(r).$$
Here $\mathbb{U}_k, \mathbb{V}_k \in o(r^k)$ satisfy
$$-\frac{\mathbb{U}_k''}{2} + \frac{\mathbb{V}_k}{r} \mathbb{U}_k - \frac{z}{r} \mathbb{U}_k - \omega \mathbb{U}_k = c_k r^{k-1} + O(r^k), r > 0$$
$$(6.10) \qquad -\mathbb{V}_k'' - 4\pi \frac{\mathbb{U}_k^2}{r} = \tilde{c}_k r^{k-1} + O(r^k), \ r > 0,$$

and these relations imply

$$\mathbb{U}_k(r) = \frac{c_k r^{k+1}}{k(k+1)} + O(r^{k+2}), \mathbb{V}_k(r) = \frac{\tilde{c}_k r^{k+1}}{k(k+1)} + O(r^{k+2})$$

near $r = 0$ and these relations imply

$$\mathbb{U}_k(r), \mathbb{V}_k(r) \in C^{k+1}([0, 1]).$$

□

Our next step is to obtain the decay of the solution. We look for soliton type solutions $u$ to (1.7), i.e. very regular solutions decaying rapidly at infinity. Our next step is to obtain a very rapid decay of the radial field $u(|x|)$ at infinity.

**Lemma 6.5.** *If the assumption* (6.1) *is satisfied, then*

(6.11) $$\mathbb{U} \in H^k((1, +\infty)), \mathbb{V}' \in H^{k-1}((1, +\infty)),$$

*and*

(6.12) $$|\mathbb{U}'(r)|^2 + |\mathbb{U}(r)|^2 \leq \frac{C}{r^k}, \qquad 0 \leq \mathbb{V}'(r) \leq \frac{C}{r^k}$$

*for each integer $k \geq 2, r \geq 1$.*

PROOF. The Sobolev embedding and Lemma 6.2 imply that

(6.13) $$\int_0^{+\infty} |\mathbb{U}(r)|^2 dr + \int_0^{+\infty} |\mathbb{U}'(r)|^2 dr \leq C \|u\|^2_{H^1(\mathbb{R}^3)},$$
$$\int_0^{+\infty} |\mathbb{V}'(r)|^2 dr \leq C \|\varphi\|^2_{\mathcal{D}^{1,2}(\mathbb{R}^3)}.$$

Note that we have used the Hardy inequality (see Theorem 330 in [14] or Remark 1, Section 3.2.6 in [21])

(6.14) $$\int_0^{+\infty} |f(r)|^2 dr \leq C \int_0^{+\infty} |f'(r)|^2 \, r^2 dr$$

in the above estimates.

Hence

$$\mathbb{U} \in H^1((0, +\infty)), \mathbb{V}' \in L^2((0, +\infty)).$$

Proceeding further inductively we find (6.11).

The above properties and the Sobolev embedding imply

(6.15) $$\lim_{r \to +\infty} |\mathbb{U}(r)| = 0, \lim_{r \to +\infty} |\mathbb{U}'(r)| = 0,$$

In a similar way we get

(6.16) $$\lim_{r \to +\infty} |\mathbb{V}'(r)| = 0,$$

We can improve the last property. Indeed, integrating the second equality in (6.8) we find

(6.17) $$\mathbb{V}'(r) = \int_r^\infty \frac{\mathbb{U}^2(\tau)}{\tau} d\tau.$$

Since

(6.18) $$\int_r^\infty \mathbb{U}^2(\tau) d\tau \leq C,$$

we get

$$0 \leq \mathbb{V}'(r) \leq \frac{C}{r}. \tag{6.19}$$

Our next step is to obtain weighted Sobolev estimates. From the first equation in (6.8) we have

$$\frac{\mathbb{U}''}{2}(r) + \omega \mathbb{U}(r) = \mathbb{F}(r),$$

$$\mathbb{F}(r) = \frac{\mathbb{V}}{r}\mathbb{U} - \frac{z}{r}\mathbb{U} \tag{6.20}$$

Since the initial data for $\mathbb{U}$ are

$$\mathbb{U}(0) = 0, \ \mathbb{U}'(0) = a_1, \tag{6.21}$$

we have the following integral equation satisfied by $\mathbb{U}$

$$\mathbb{U}(r) = \sinh(\sqrt{-2\omega}r)a_1 + \int_0^r \sinh(\sqrt{-2\omega}(r-\rho))\mathbb{F}(\rho)d\rho. \tag{6.22}$$

It is easy to see that for $= z$ the function $\mathbb{F}$ satisfies the estimate

$$\mathbb{F}(r) = O(r^{-2}), \ r \geq 1. \tag{6.23}$$

Then the condition (6.15) and simple qualitative study of the integral equation in (6.22) guarantees that

$$a_1 + \int_0^\infty e^{\sqrt{-2\omega}\rho}\mathbb{F}(\rho)d\rho = 0.$$

This fact enables one to represent $\mathbb{U}$ as follows

$$\mathbb{U}(r) = e^{-\sqrt{-2\omega}r}a_1 - \tag{6.24}$$

$$- \int_r^\infty e^{\sqrt{-2\omega}(r-\rho)}\mathbb{F}(\rho)d\rho - \int_0^r e^{-\sqrt{-2\omega}(r-\rho)}\mathbb{F}(\rho)d\rho.$$

Since

$$\mathbb{F}(r) = \frac{\mathbb{V}}{r}\mathbb{U} - \frac{z}{r}\mathbb{U} = O(r^{-2})$$

we obtain

$$|\mathbb{U}(r)| = O(r^{-1}).$$

This estimate implies a stronger version of (6.18)

$$\int_r^\infty \mathbb{U}^2(\tau)d\tau \leq \frac{C}{r}, \tag{6.25}$$

and from (6.17) we improve (6.19) as follows

$$0 \leq \mathbb{V}'(r) \leq \frac{C}{r^2}. \tag{6.26}$$

This argument shows that combining (6.19) and (6.20) we can obtain inductively

$$\sum_{j=0}^k |\mathbb{U}^{(j)}(r)|^2 \leq \frac{C}{r^n} \tag{6.27}$$

and

$$\sum_{j=1}^k |\mathbb{V}^{(j)}(r)|^2 \leq \frac{C}{r^n} \tag{6.28}$$

for any integers $k \geq 1$ and $n \geq 2$. $\square$

PROOF OF THEOREM 1.4. The proof is immediate consequence of Lemmas 6.5 and 6.2, and change of variables (6.7). $\square$

## 7. PROOF OF THEOREM 1.5

In this section we shall prove Theorem 1.5 We begin proving some lemmas. Define the functional

$$I(u,\omega) \doteq \frac{1}{4}\int_{\mathbb{R}^3}|\nabla u|^2 dx + \pi \int_{\mathbb{R}^3}|\nabla \Delta^{-1}u^2|^2 dx-$$
$$(7.1) \qquad -\frac{1}{2}\int_{\mathbb{R}^3}V(x)|u|^2 dx - \frac{\omega}{2}\int_{\mathbb{R}^3}|u|^2 dx,$$

for each $(u,\omega) \in H_r^1(\mathbb{R}^3) \times \mathbb{R}$. There results

$$\frac{\partial I}{\partial u}(u,\omega) = -\frac{1}{2}\Delta u - 4\pi(\Delta^{-1}u^2)u - V(x)u - \omega u,$$
$$\frac{\partial I}{\partial \omega}(u,\omega) = -\frac{1}{2}\int_{\mathbb{R}^3}|u|^2 dx,$$
$$\frac{\partial^2 I}{\partial u^2}(u,\omega)h = -\frac{1}{2}\Delta h - 4\pi(\Delta^{-1}u^2)h - 8\pi\Delta^{-1}(hu)u - V(x)h - \omega h, \quad h \in H_r^1(\mathbb{R}^3),$$
$$\frac{\partial^2 I}{\partial u \partial \omega}(u,\omega) = -u,$$
$$\frac{\partial^2 I}{\partial \omega^2}(u,\omega) = 0.$$

Let

$$\nabla I : H_r^1(\mathbb{R}^3) \times \mathbb{R} \longrightarrow \left(H_r^1(\mathbb{R}^3)\right)' \times \mathbb{R}, \qquad \nabla I(u,\omega) = \begin{pmatrix} \dfrac{\partial I}{\partial u}(u,\omega) \\ \dfrac{\partial I}{\partial \omega}(u,\omega) \end{pmatrix}$$

be the Jacobian matrix of $I$ and

$$HI(u,\omega) : H_r^1(\mathbb{R}^3) \times \mathbb{R} \to \left(H_r^1(\mathbb{R}^3)\right)' \times \mathbb{R}, HI(u,\omega) = \begin{pmatrix} \dfrac{\partial^2 I}{\partial u^2}(u,\omega) & \dfrac{\partial^2 I}{\partial u \partial \omega}(u,\omega) \\ \dfrac{\partial^2 I}{\partial u \partial \omega}(u,\omega) & \dfrac{\partial^2 I}{\partial \omega^2}(u,\omega) \end{pmatrix}$$

be the Hessian matrix of $I$ in $(u,\omega)$. More precisely

$$HI(u,\omega)(h,k) = \begin{pmatrix} \dfrac{\partial^2 I}{\partial u^2}(u,\omega)h + \dfrac{\partial^2 I}{\partial u \partial \omega}(u,\omega)k \\ \dfrac{\partial^2 I}{\partial u \partial \omega}(u,\omega)h + \dfrac{\partial^2 I}{\partial \omega^2}(u,\omega)k \end{pmatrix} =$$

$$(7.2) \qquad = \begin{pmatrix} -\dfrac{1}{2}\Delta h - 4\pi(\Delta^{-1}u^2)h - 8\pi(\Delta^{-1}(hu))u - V(x)h - \omega h - ku \\ -\displaystyle\int_{\mathbb{R}^3} uh\, dx \end{pmatrix},$$

for each $u$, $h \in H^1_r(\mathbb{R}^3)$ and $k$, $\omega \in \mathbb{R}$. Finally denote

(7.3) $$B' \doteq B \cap H^1_r(\mathbb{R}^3).$$

**Lemma 7.1.** *Let $u_0 \in B'$ (see (7.3)) be a critical point of $J\big|_{B'}$ that corresponds to the minimum*

(7.4) $$\omega_0 = \inf_{u \in H^1 \setminus \{0\}, \|u\|^2_{L^2} = N} J(u),$$

*namely*

$$0 = J\big|'_{B'}(u_0) = J\big|'_{B}(u_0) = J'(u_0) - \omega_0 u_0.$$

*The operator*

$$h \in \left\{ h \in H^1_r(\mathbb{R}^3), \int h(x) u_0(x) dx = 0 \right\} \longmapsto \frac{\partial^2 I}{\partial u^2}(u_0, \omega_0) h \in \left(H^1_r(\mathbb{R}^3)\right)'$$

*has a trivial kernel and*

(7.5) $$\left\langle \frac{\partial^2 I}{\partial u^2}(u_0, \omega_0) h \big| h \right\rangle = 0, \quad \int h(x) u_0(x) dx = 0 \implies h \equiv 0.$$

PROOF. Repeating the qualitative argument in the proof of Lemma 6.5, we see that any solution of

$$\frac{\partial^2 I}{\partial u^2}(u_0, \omega_0) h = 0$$

decays rapidly at infinity and it is smooth as a function of $r \geq 0$. Another interpretation of the first eigenvalue $\omega_0 = \omega(N) < 0$ is the following one

(7.6) $$\omega_0 = N \inf_{u \in H^1 \setminus \{0\}} \frac{J(u)}{\|u\|^2_{L^2}}.$$

Let

$$\left\langle \frac{\partial^2 I}{\partial u^2}(u_0, \omega_0) h \big| h \right\rangle = 0$$

for some $h \in H^1$ orthogonal (in $L^2$) to $u_0$. Take $u_0 + i\varepsilon h$ with $\varepsilon > 0$ small enough (will be chosen later on). Then a simple calculation implies

$$\frac{J(u_0 + i\varepsilon h)}{\|u_0 + i\varepsilon h\|^2_{L^2}} = \frac{J(u_0) + o(\varepsilon^2)}{\|u_0\|^2 - \varepsilon^2 \|h\|^2} = \frac{J(u_0)}{N} + \varepsilon^2 \frac{\|h\|^2}{N} J(u_0) + o(\varepsilon^2).$$

Hence, the assumption $\|h\| \neq 0$ will contradict the fact that $\omega_0$ is defined as the minimum in (7.6). This completes the proof. □

**Lemma 7.2.** *Let $u_0 \in B'$ (see (7.3)) be a critical point of $J\big|_{B'}$ that corresponds to the minimum as in the previous Lemma. The operator*

$$(h, k) \in H^1_r(\mathbb{R}^3) \times \mathbb{R} \longmapsto HI(u_0, \omega_0)(h, k) \in \left(H^1_r(\mathbb{R}^3)\right)' \times \mathbb{R}$$

*is invertible.*

PROOF. Let $u_0 \in B'$ be a critical point of $J\big|_{B'}$ with multiplier $\omega_0$ as in the previous lemma, call

$$A \doteq \frac{\partial^2 I}{\partial u^2}(u_0, \omega_0).$$

We begin proving that $HI(u_0, \omega_0)$ is injective. Let $h \in H_r^1(\mathbb{R}^3)$ and $k \in \mathbb{R}$ such that

(7.7) $$HI(u_0, \omega_0)(h, k) = 0,$$

we have to prove that

(7.8) $$h = k = 0.$$

By (7.2) and (7.7), we have

(7.9) $$Ah - ku_0 = 0, \quad -\int_{\mathbb{R}^3} u_0 h \, dx = 0.$$

Multiplying the first of (7.9) by $h$ and integrating on $\mathbb{R}^3$, we have

$$\int_{\mathbb{R}^3} (Ah) h \, dx = -k \int_{\mathbb{R}^3} u_0 h \, dx = 0,$$

and by (7.5) and the definition of $A$

(7.10) $$h \equiv 0.$$

On the other hand, multiplying the first of (7.9) by $u_0$ and integrating on $\mathbb{R}^3$, since $u_0 \in B'$, we have

(7.11) $$kN = k\int_{\mathbb{R}^3} u_0^2 \, dx = \int_{\mathbb{R}^3} (Ah) u_0 \, dx = 0.$$

Since (7.8) is direct consequence of (7.10) and (7.11), $HI(u_0, \omega_0)$ is injective.

We prove that $HI(u_0, \omega_0)$ is surjective. Observe that the operator $A$ is selfadjoint, indeed

(7.12)
$$(Ah, f)_{L^2} = \frac{1}{2}\int_{\mathbb{R}^3}(\nabla h, \nabla f) \, dx - 4\pi \int_{\mathbb{R}^3}(\Delta^{-1}u^2) h f \, dx+$$

(7.13) $$+ 8\pi \int_{\mathbb{R}^3} \left(\nabla \Delta^{-1}(hu), \nabla \Delta^{-1}(fu)\right) dx - \int_{\mathbb{R}^3} V(x) h f \, dx - \omega \int_{\mathbb{R}^3} h f \, dx,$$

(7.14)

for each $h, f \in H_r^1(\mathbb{R}^3)$. Moreover, also the operator $HI(u_0, \omega_0)$ is selfadjont, indeed

(7.15) $$\left(HI(u_0, \omega_0)(h, k), (f, \alpha)\right)_{L^2 \times \mathbb{R}} = \left(\left(Ah - ku_0, -(u_0, h)_{L^2}\right), (f, \alpha)\right)_{L^2 \times \mathbb{R}} =$$

(7.16) $$= \left(Ah - ku_0, f\right)_{L^2} - \alpha(u_0, h)_{L^2} =$$

(7.17) $$= \left(Ah, f\right)_{L^2} - k\left(u_0, f\right)_{L^2} - \alpha(u_0, h)_{L^2}$$

and

(7.18) $$\left(HI(u_0, \omega_0)(f, \alpha), (h, k)\right)_{L^2 \times \mathbb{R}} = \left(\left(Af - \alpha u_0, -(u_0, f)_{L^2}\right), (h, k)\right)_{L^2 \times \mathbb{R}} =$$

(7.19) $$= \left(Af - \alpha u_0, h\right)_{L^2} - k(u_0, f)_{L^2} =$$

(7.20) $$= \left(Af, h\right)_{L^2} - \alpha(u_0, h)_{L^2} - k\left(u_0, f\right)_{L^2},$$

since $A$ is selfadjoint
$$\big(HI(u_0,\,\omega_0)(h,k),(f,\alpha)\big)_{L^2\times\mathbb{R}} = \big(HI(u_0,\,\omega_0)(f,\alpha),(h,k)\big)_{L^2\times\mathbb{R}}$$
for each $h,\,f\,\in H_r^1(\mathbb{R}^3)$ and $k,\,\alpha\,\in\mathbb{R}$. Since $HI(u_0,\omega_0)$ is injective and selfadjoint, there results

(7.21) $\quad\mathrm{Im}\big(HI(u_0,\,\omega_0)\big) = \Big(\mathrm{Ker}\big(HI(u_0,\,\omega_0)^*\big)\Big)^\perp =$

(7.22) $\quad\qquad\qquad = \Big(\mathrm{Ker}\big(HI(u_0,\,\omega_0)\big)\Big)^\perp = H_r^1(\mathbb{R}^3)\times\mathbb{R},$

then $HI(u_0,\omega_0)$ is surjective. The claim is direct consequence of the Closed Graph Theorem. $\square$

**Lemma 7.3.** *The critical points of the functional $J\big|_{B'}$ that correspond to the minimum are isolated, i.e. for each $u\in B'$ critical point of $J\big|_{B'}$, with the Lagrange multiplier satisfying (7.4), there exists a neighborhood $U\subset H^1(\mathbb{R}^3)$ of $u$ such that any element of $B'\cap U$ is not a critical point of it.*

PROOF. Let $u_0\in B'$ be a critical point of $J\big|_{B'}$ corresponding to the minimum as in the previous lemmas, then
$$0 = J\big|_{B'}'(u_0) = J\big|_B'(u_0) = J'(u_0) - \omega_0 u_0 = \frac{\partial I}{\partial u}(u_0,\,\omega_0)$$
and since $u_0\in B'$
$$\frac{\partial I}{\partial\omega}(u_0,\,\omega_0) = -\frac{1}{2}\int_{\mathbb{R}^3}u_0^2 dx = -\frac{N}{2},$$
we have
$$\nabla I(u_0,\,\omega_0) = \begin{pmatrix}0\\-N/2\end{pmatrix}.$$
By Lemma 7.2 and the Implicit Function Theorem there exist $U\subset H_r^1(\mathbb{R}^3)$ neighborhood of $u_0$, $\Omega\subset\mathbb{R}$ neighborhood of $\omega_0$, $W\subset\big(H_r^1(\mathbb{R}^3)\big)'\times\mathbb{R}$ neighborhood of $\left(0,\,-\dfrac{N}{2}\right)$ and $G:W\longrightarrow U\times\Omega$ such that

$$G\big(\nabla I(u,\,\omega)\big) = (u,\,\omega), \qquad\qquad (u,\omega)\in U\times\Omega,$$
(7.23) $\quad\nabla I\big(G(f,\,\alpha)\big) = (f,\,\alpha), \qquad\qquad (f,\,\alpha)\in W.$

Assume, by absurd, that $u_0$ is not isolate, namely there exists a sequence $\{u_k\}\subset B'$ of critical points of $J\big|_{B'}$, such that

(7.24) $\qquad\qquad u_k\neq u_0, \qquad u_k\longrightarrow u_0 \quad\text{in } H^1(\mathbb{R}^3).$

Moreover, there exists a sequence $\{\omega_k\}\subset\mathbb{R}$ such that
$$0 = J\big|_{B'}'(u_k) = J'(u_k) - \omega_k u_k = \frac{\partial I}{\partial u}(u_k,\,\omega_k).$$
Since $u_k\in B'$ and by (7.24), we have

(7.25) $\qquad\qquad \omega_k = \big\langle J'(u_k)\big|u_k\big\rangle \longrightarrow \big\langle J'(u_0)\big|u_0\big\rangle = \omega_0.$

By (7.24) and (7.25), there exists $k_0\in\mathbb{N}$ such that
$$(u_k,\,\omega_k)\in U\times\Omega, \qquad k\geq k_0.$$

Finally, fixed $k \geq k_0$, since

$$\nabla I(u_k, \omega_k) = \begin{pmatrix} 0 \\ -N/2 \end{pmatrix},$$

by (7.23), we have

$$(u_k, \omega_k) = G(\nabla I(u_k, \omega_k)) = G\begin{pmatrix} 0 \\ -N/2 \end{pmatrix} = G(\nabla I(u_0, \omega_0)) = (u_0, \omega_0).$$

Since this contradicts (7.24), the claim is done. □

**Lemma 7.4.** *The first eigenvalue of the operator $J\big|'_{B'}$ (see (7.4)) is isolated, i.e. there exists a neighborhood $\Omega \subset \mathbb{R}$ of $\omega_0$ such that any element of $\Omega$ is not an eigenvalue of the previous operator.*

PROOF. Assume, by absurd, that the first eigenvalue $\omega_0$ is not isolated, namely there exists a sequence $\{\omega_k\} \subset \mathbb{R}$ of eigenvalues such that

(7.26) $$\omega_k \longrightarrow \omega_0.$$

By definition, there exists $\{u_k\} \subset B'$ such that

(7.27) $$0 = J\big|'_{B'}(u_k) = J'(u_k) - \omega_k u_k, \qquad k \in \mathbb{N}.$$

Observe that, by Lemma 4.3, $\omega_k, \omega_0 < 0$, then there exists $\varepsilon > 0$ such that

(7.28) $$\omega_k, \omega_0 \leq -\varepsilon, \qquad k \in \mathbb{N}.$$

Moreover, by Lemma 3.3 and since $\{u_k\} \subset B'$

(7.29) $$-\infty < \min_{u \in H^1_r(\mathbb{R}^3)} J(u) \leq J(u_k) \leq \sup_k \frac{\omega_k}{2} \leq -\frac{\varepsilon}{2},$$

then $\{J(u_k)\}$ is bounded and, by (7.27),

(7.30) $$J\big|'_{B'}(u_k) \longrightarrow 0.$$

By the Palais-Smale Condition (see Lemma 5.1) there exists $u_0 \in B'$ such that, passing to a subsequence,

$$u_k \longrightarrow u_0, \qquad \text{in } H^1(\mathbb{R}^3).$$

By (7.26) and (7.27),

$$0 = J\big|'_{B'}(u_0) = J'(u_0) - \omega_0 u_0,$$

namely $u_0$ is a not isolated critical point of the functional $J\big|_{B'}$. Since this contradicts Lemma 7.3 the proof is done. □

PROOF OF THEOREM 1.5. Since

$$F(u, 4\pi\Delta^{-1}u^2) = J(u)$$

for all $u \in H^1(\mathbb{R}^3)$, by Lemmas 7.3 and 7.4 the claim is done. □

## 8. Appendix

Here we shall prove for completeness the relation (4.10). First, for the partial case of space dimensions $n = 3$ we need the following relation ( a generalization of this relation for space dimensions $n \geq 3$ can be found in [1]).

**Lemma 8.1.** ( see [1]) *If $f(x) = f(|x|)$ is an $L^\infty(\mathbb{R}^3)$ function, then for any $r > 0$ and $x \neq 0$ we have the relation*

$$(8.1) \qquad \int_{\mathbb{S}^2} f(|x + r\omega|) d\omega = \frac{2\pi}{|x|r} \int_{||x|-r|}^{|x|+r} f(\lambda) \lambda d\lambda.$$

PROOF . It is sufficient to consider only the case $x = (0, 0, |x|)$ and to pass to polar coordinates

$$\omega_1 = \sin\theta \cos\varphi \ , \ \omega_2 = \sin\theta \sin\varphi \ , \ \omega_3 = \cos\theta.$$

Then $d\omega = \sin\theta\, d\theta\, d\varphi$ and

$$\int_{\mathbb{S}^2} f(|x+r\omega|) d\omega = 2\pi \int_0^\pi f\left(\sqrt{|x|^2 + r^2 + 2|x|r\cos\theta}\right) \sin\theta\, d\theta.$$

Making the change of variable

$$\theta \to \lambda = \sqrt{|x|^2 + r^2 + 2|x|r\cos\theta},$$

we complete the proof. $\square$

Now we are ready to verify (4.10).

**Lemma 8.2.** *If $v(x) = v(|x|)$ is a radial $C_0^\infty(\mathbb{R}^3)$ function, then the solution of the equation*

$$\Delta u = v$$

*can be represented as follows*

$$(8.2) \qquad 4\pi u(x) = -\int_{\mathbb{R}^3} v(|y|) \frac{dy}{\max(|x|,|y|)}, \qquad x \in \mathbb{R}^3.$$

PROOF . Starting with the classical representation

$$4\pi u(x) = \int_{\mathbb{R}^3} |x-y|^{-1} v(|y|) dy,$$

we introduce polar coordinates $r = |y|, \omega = y/|y|$ apply Lemma 8.1 and find

$$u(x) = -\frac{1}{2|x|} \int_0^\infty \left(\int_{||x|-r|}^{|x|+r} d\lambda\right) v(r) r dr.$$

Note that the right side of (8.2) becomes

$$-4\pi \int_0^\infty v(r) \frac{r^2 dr}{\max(|x|,r)}.$$

Using the fact that

$$\frac{1}{|x|r} \int_{||x|-r|}^{|x|+r} d\lambda = \frac{2}{\max(|x|,r)},$$

we obtain (8.2) and this completes the proof. $\square$

Using the relation

$$4\pi u(x) = -\int_0^r v(\rho) \frac{\rho^2 d\rho}{r} - \int_r^\infty v(\rho) \rho d\rho, \quad r = |x|$$

and differentiating with respect to $r = |x|$, we arrive at the following.

**Lemma 8.3.** *If $v(x) = v(|x|)$ is a radial $C_0^\infty(\mathbb{R}^3)$ function, then the solution of the equation*

$$\Delta u = v$$

*satisfies the relation*

(8.3) $$4\pi \frac{\partial \Delta^{-1} v}{\partial r}(x) = \int_{|y|<r} \frac{v(y)}{|x|^2} dy,$$

*for each $x \in \mathbb{R}^3$, $x \neq 0$.*

## Acknowledgments

The authors would like to thank Dott. Simone Secchi for for stimulating conversations.

Giuseppe Maria Coclite, S.I.S.S.A., Via Beirut 2 - 4, 34013 Trieste, Italy, E-mail: coclite@sissa.it

Vladimir Georgiev, Dipartimento di Matematica, Università degli Studi di Pisa, Via F. Buonarroti 2, 56100 Pisa, Italy, E-mail: georgiev@dm.unipi.it